\documentclass[11pt]{article}

\usepackage{fullpage,amssymb,amsbsy}

\newcommand{\norm}[1]{\| #1\|}
\newtheorem{theo}{Theorem}
\newtheorem{lemma}[theo]{Lemma}
\newtheorem{cor}[theo]{Corollary}
\newtheorem{prop}[theo]{Proposition}
\newtheorem{defn}[theo]{Definition}
\newenvironment{Proof}{\removelastskip \vskip12pt plus 1pt \noindent
{\em Proof.\/}\rm }{\hfill$\Box$ \vskip12pt plus 1pt} 

\newcommand{\R}{{\mathbb R}}

\title{Minimal Travelling Wave Speed and Explicit Solutions in
  Monostable Reaction-Diffusion Equations}

\author{E.C.M. Crooks, Department of Mathematics, College of
  Science,\\ Swansea University, Bay Campus, Swansea SA1 8EN
  \texttt{e.c.m.crooks@swansea.ac.uk}\\ and\\
  M. Grinfeld, Department of Mathematics and Statistics,\\
  University of Strathclyde, Glasgow G1 1XH
  \texttt{m.grinfeld@strath.ac.uk}}

\begin{document}

\maketitle


\begin{abstract}
  We investigate the connection between the existence of an explicit
  travelling wave solution and the travelling wave with minimal speed
  in a scalar monostable reaction-diffusion equation.
\end{abstract}

\section{Introduction}

In this short paper we investigate the somewhat puzzling connection
between the existence of an explicit travelling wave solution and the
travelling wave with minimal speed in a monostable reaction-diffusion
equation. More precisely, it often happens that the explicitly
computable travelling wave solution is the solution with minimal
speed. Moreover, for parameter-dependent problems with a
parameter-dependent family of explicit solutions, it is common for
there to be in fact a switching between the minimal speed being given
by the explicit solution for some parameters, while for others it is
given by the so-called linear speed, defined as the minimal value for
which the problem linearised about the unstable steady state has a
suitable eigenvalue.  For a particular set of equations, of a type
often encountered in applications, we formulate sufficient conditions
for each of these phenomena to occur.

The plan of the paper is as follows. In this section, we introduce
scalar monostable reaction-diffusion equations, define what we mean by
a minimal speed, and discuss the linear (pulled) and the non-linear
(pushed) regimes. 

In section \ref{sol}, we define the set of exactly solvable
equations and prove a result connecting the minimal wave speed and
the speed of the explicit travelling wave solution.

Finally, in section \ref{ex} we consider conditions for the
exchange of minimality between the linear minimal speed and the speed
of the explicit travelling wave solution. 

Throughout, in all our proofs we only use two tools: the variational
principle due to Hadeler and Rothe \cite{HR} and the integrability
characterisations of the minimal speed proved by Lucia, Muratov and Novaga in \cite{LMN}. 

We consider reaction-diffusion equations of the form
\begin{equation}\label{rd}
u_t= u_{xx}+f(u,\beta),
\end{equation}
where $\beta \in \R$ is a parameter, and $f$ is a 
monostable nonlinearity, {\em {\em i.e.}, }
\[
f(0,\beta)=f(1,\beta)=0, \; \; f'(0,\beta)>0, \;\; f'(1, \beta)<0, \;\; f(u,\beta)>0 
\hbox{ for } u \in (0,1).
\]

In the travelling wave frame $z=x-ct$, $c \geq 0$, setting $U(z)=u(x,t)$,
  and denoting derivatives with respect to  $z$ by primes, (\ref{rd}) becomes
\begin{equation}\label{tw}
    -c U' = U'' + f(U).
\end{equation}    

We seek {\bf monotone fronts} connecting $1$ and $0$,
{\em i.e.},  solutions $U(z)$ of (\ref{tw}) such that
\[
  \lim_{z\rightarrow -\infty} U(z)=1 \hbox{ and  } \lim_{z\rightarrow
    \infty} U(z)=0.
\]

Linearisation around the rest point with $U=0$ shows that there cannot
be any monotone fronts connecting $1$ and $0$ for $c < c_l :=
2\sqrt{f'(0)}$. Phase plane analysis shows that there exists $c_{min}
\geq c_l$ such that there exists a monotone front for all $c \geq
c_{min} \geq c_l$. Determining $c_{min}$ is often of interest in
applications, see e.g. \cite{CGM} for a discussion.

\begin{defn} If $c_{min}=c_l$, we say that we are in
  the case of {\em linear selection mechanism\/} (``pulled case'') and
  if $c_{min} > c_l$, of {\em nonlinear selection mechanism} (``pushed
  case'').
\end{defn}

The basis of almost all analysis of monotone fronts in the scalar
monostable case (\ref{tw}) is the following construction: As $U(z)$ is
a monotone solution, its derivative is a well-defined function of
$U$. Set $F(U):= -U'$. Note that $F(U)$ is non-negative. Also,
$F(0)=F(1)=0$. Now,
\[
  F(U)'=(-U')'=-U''.
\]
On the other hand, by the chain rule,
\[
  F(U)' = \frac{dF}{dU}U'= - \frac{dF}{dU}F.
\]
Hence the problem of solving $U''+cU'+f(U)=0$ with the conditions that
$ \lim_{z\rightarrow -\infty} U(z)=1$ and $\lim_{z\rightarrow \infty}
U(z)=0$ is equivalent to solving 
\begin{equation}\label{feq}
 F \frac{dF}{dU} -cF + f(U) =0, \; \; F(0)=F(1)=0. 
\end{equation}

Using this construction, we have the Hadeler--Rothe variational
principle \cite{HR}:

\begin{equation}\label{HR}
  c_{min} = \inf_{g \in {\cal G}} \sup_{0<U<1} \left\{ g'(U) +
    \frac{f(U)}{g(U)} \right\}, 
\end{equation}
where 
\begin{equation}\label{setG}
  {\cal G} = \{ g \in C^1([0,1]) \, | \, g(U) >0 \hbox{ for } 0<U<1,
  \, g(0)=0, \; g'(0)>0 \}.
\end{equation}

\section{Exact solvability}\label{sol}

We are interested in the situation when (\ref{tw}) has a solution
$U(z)$ that can be determined by quadratures. A sufficient
condition is:

\begin{lemma} \label{solv}
The travelling wave equation of (\ref{feq})
  is solvable by quadratures if $f$ can be written in the form
\[
f(u,\beta)= h(u) \bigg( A(\beta)-B(\beta)h'(u) \bigg),
\]
where $h(0)=h(1)=0$, $h(u)\geq 0$, (without loss of generality $h'(0)=1$),
$A(\beta)>$, $B(\beta)>0$, and  for all $u \in [0,1]$,
$A(\beta)-B(\beta)h'(u)>0$.
\end{lemma}

\begin{Proof} In this case a solution of (\ref{feq}) is $F(U)=\gamma h(U)$ with 
\begin{equation}\label{gamma}
\gamma= \sqrt{B(\beta)},
\end{equation}
from which $U$ can be computed by quadratures.
\end{Proof}

In fact, we can compute the speed without solving for the front profile: 
\begin{equation}\label{cnl}
c := c_{nl} (\beta)= \frac{A(\beta)}{\sqrt{B(\beta)}}.
\end{equation}

We will describe as  the {\bf solvable
  case} the situation in which the nonlinearity $f(u,\beta)$
satisfies the conditions of Lemma \ref{solv}. In the solvable case, we have that 
\begin{equation}\label{cls}
c_l = 2 \sqrt{A(\beta)-B(\beta)}.
\end{equation}
Note that the fact that  $A(\beta)>B(\beta)$ follows from the conditions of Lemma
\ref{solv}. 

Of course, by the definition of minimal speed, we always have that
\begin{equation}\label{apr}
  c_{min} (\beta) \leq c_{nl}(\beta) = \frac{A(\beta)}{\sqrt{B(\beta)}}.
\end{equation}

\section{Minimality Exchange}\label{ex}

In this section, for a nonlinearity $f(u, \beta)$ of solvable type, we
investigate conditions under which there exists a value $\beta^*$, 
such that for values $\beta$ to one side of $\beta^*$,
$c_{min}(\beta)=c_l(\beta)$, and for values of $\beta$ to the other
side of $\beta^*$, $c_{min}(\beta)=c_{nl}(\beta)$, so that at
$\beta^*$ minimality is exchanged between $c_l(\beta)$ and
$c_{nl}(\beta)$. This is what we call a {\bf minimality exchange}. Examples, two of which we outline below, 
are discussed in \cite{HR,LMN} and the isotropic case of \cite{CGM},
which is also investigated in \cite{GK, SHH}.

First note that   for a minimality exchange, $c_l(\beta)$
and $c_{nl}(\beta)$ must clearly intersect. Therefore the
equation
\[
2 \sqrt{A(\beta)-B(\beta)} = \frac{A(\beta)}{\sqrt{B(\beta)}}
\]
must have a solution, which is equivalent to demanding the
existence of $\beta^*$ such that $A(\beta^*)=2B(\beta^*)$. 

Hence, for instance, in any equation $(\ref{rd})$ with solvable
$f(u,\beta)$ such that $A(\beta)=2B(\beta)+1$, there can never be a
minimality exchange between the linear and the nonlinear speeds. 

Before continuing with the analysis, we present two  concrete examples of minimality exchange. 
In \cite[Eq. (27)]{HR}, Hadeler and Rothe consider the nonlinearity
\[
  f(u, \beta) = u(1-u)(1+\beta u), \; \beta \geq -1, 
\]
which can be put into the framework of Lemma \ref{solv} by setting $h(u)=u(1-u)$, so that
\[
  f(u,\beta) = h(u)(A(\beta)-B(\beta) h'(u)),
\]
where
\[
  A(\beta) = 1 +\frac{\beta}{2}, \; B(\beta) = \frac{\beta}{2}. 
\]
The solution of $A(\beta)=2B(\beta)$ is therefore $\beta^* =2$, the
nonlinear speed is
\[
  c_{nl}(\beta)=\frac{2+\beta}{\sqrt{2\beta}}, 
\]    
and it is shown in \cite{HR} that a minimality exchange
  occurs at $\beta = \beta^*$,  with $c_{min}(\beta)=c_l(\beta)$ for $\beta < \beta^*$ and $c_{min}(\beta)=c_{nl}(\beta)$
  for $\beta > \beta^*$. 

Our second example is given by the isotropic case of \cite{CGM}, where 
\[ f(u, \beta) =\frac{\sin (\pi u)}{2\pi} \left[1 -\beta \cos(\pi
  u)\right],
  \]
  which fits into the framework of Lemma \ref{solv} by setting $\displaystyle{h(u)=\frac{\sin(\pi u)}{\pi}}$, so that
\[
  f(u,\beta) = h(u)(A(\beta)-B(\beta) h'(u)),
\]
where
\[
  A(\beta) = \frac{1}{2} \; B(\beta) = \frac{\beta}{2}.
  \]
The equation  $A(\beta)=2B(\beta)$ then has solution  $\beta^*=\frac{1}{2}$, the nonlinear speed is
  \[ c_{nl}(\beta) = \frac{1}{\sqrt{2 \beta}},\]
  and it is proved in \cite{CGM, GK} that here too,  a minimality exchange
  occurs at $\beta = \beta^*$,  again with $c_{min}(\beta)=c_l(\beta)$ for $\beta < \beta^*$ and $c_{min}(\beta)=c_{nl}(\beta)$
  for $\beta > \beta^*$. 
  

\bigskip

We now establish our general results, starting with a sufficient condition for nonlinear selection. 

\begin{lemma}\label{lnl}
  For all $\beta$ such that $A(\beta) < 2 B(\beta)$,
  $c_{min}(\beta)=c_{nl}(\beta)$.
\end{lemma}

\begin{Proof}
For any $c>0$, denote by $H^1_c(\R)$ the completion of $C_0^\infty(\R)$
with respect to the norm
\[
\norm{u}_{1,c} = \norm{u}_c + \norm{u_x}_c, \hbox{ where }
\norm{u}_c^2 = \int_\R e^{cx} u^2(x) \, dx.
\]
If $U(z)$ is an explicit travelling front with $-U'=F(U)=\gamma h(U)$,  by L'H\^opital's rule
\[
\lim_{z \to \infty}\frac{U'(z)}{U(z)} =\lim_{z \to \infty} - \gamma \, 
\frac{h(U(z))}{U(z)} = - \gamma\, h'(0)= - \gamma. 
\]
Hence for those values of the
parameter $\beta$ for which $c_{nl}(\beta)<2 \gamma$, $U \in
H^1_{c_{nl}(\beta)}(\R)$ and hence for such $\beta$, by Corollary 2.7 of \cite{LMN} (see
also Proposition 2 of \cite{CGM}), $c(\beta)$ is the
(nonlinear) minimal wave speed. The claim then follows by
(\ref{gamma}) and (\ref{cnl}).
\end{Proof}

To formulate our next results, we set 
\[
L = \max_{u \in (0,1)}h'(u) \geq 1. 
\]  

We adapt some arguments from \cite{CGM}.

\begin{prop} If $A(\beta) > 2L
B(\beta)$,
\begin{equation}
\label{bigLest}
c_{min}(\beta) \leq 2\sqrt{L}\sqrt{A(\beta)-LB(\beta)},
\end{equation}
and in particular, 
$$c_{min}(\beta) \neq c_{nl}(\beta).$$
\end{prop} 

\begin{Proof} Recall from  Hadeler and Rothe \cite{HR} (see also \cite{CGM}, equation (11)) that
\begin{equation}
\label{minmax}
c_{min}(\beta) = \inf_{g \in \Lambda} \; \sup_{U \in (0,1)} \;\left\{ g'(U) + 
\frac{f(U, \beta))}{g(U)} \right\},
\end{equation}
where
\begin{equation}
\label{lambdadef}
\Lambda = \{ g \in C^1([0,1]): g(U)>0 \; \mbox{if} \; U \in (0,1), \; g(0)=0, 
\; g'(0)>0\}. 
\end{equation}

Hence taking  $g(U)=
\nu h(U)$, $\nu >0$, yields that 
\[
c_{min}(\beta)  \leq \inf_{\nu>0} \sup_{U \in (0,1)} \left\{ 
\nu h'(U) + \frac{A(\beta)}{\nu} -\frac{B(\beta)}{\nu}h'(U) \right\}. 
\]
To understand
\[
\sup_{U \in (0,1)} \left\{ 
\left( \nu -\frac{B(\beta)}{\nu}\right)  h'(U) + \frac{A(\beta)}{\nu} \right\},
\]
there are two cases:   

(i) $\nu^2 \leq B(\beta)$: Then
\[ 
\sup_{U \in (0,1)} \left\{ 
\left( \nu - \frac{B(\beta)}{\nu}\right)  h'(U) + \frac{A(\beta)}{\nu} \right\} =
\frac{A(\beta)-lB(\beta)}{\nu}+l\nu,
\]
which is monotone decreasing in $\nu$,  so \[
\inf_{\nu \leq \sqrt{B(\beta)}} \sup_{U \in (0,1)} \left\{ 
\left( \nu - \frac{B(\beta)}{\nu}\right)  h'(U) + \frac{A(\beta)}{\nu} \right\} =\frac{A(\beta)}{\sqrt{B(\beta)}}.
\]
(Note that this recovers the estimate $(\ref{apr})$ for  $c_{min}(\beta)$.)

(ii) $\nu^2 \geq B$: Then  
\[ 
\sup_{U \in (0,1)} \left\{ 
\left( \nu - \frac{B(\beta)}{\nu}\right)  h'(U) + \frac{A(\beta)}{\nu} \right\} =
\frac{A(\beta)-LB(\beta)}{\nu}+ L \nu := q(\nu).
\]
Since $A(\beta) - B(\beta) h'(u)>0$ for all $u \in [0, 1]$, it follows that $A(\beta) - LB(\beta) > 0$. 
So
differentiating $q(\nu)$ gives that its global minimum for
$\nu \in (0, \infty)$
occurs at
\[
\nu_0 := \sqrt{\frac{A(\beta)-LB(\beta)}{L}}.
\]

There are two possibilities: (a) If 
\[
\frac{A(\beta)-LB(\beta)}{L} \leq  B(\beta),
\]
the function $q(\nu)$ reaches its minimum  over $ [\sqrt{B(\beta)}, \infty)$ at the point $\nu = \sqrt{B(\beta)}$, so that
\[
\inf_{\nu  \geq \sqrt{B(\beta)}} \sup_{U \in (0,1)} \left\{ 
\left( \nu - \frac{B(\beta)}{\nu}\right)  h'(U) + \frac{A(\beta)}{\nu} 
\right\} =\frac{A(\beta)}{\sqrt{B(\beta)}},
\]
in which case we again just recover the estimate (\ref{apr}) for $c_{min}(\beta)$.

(b)  On the other hand, if  
\[
\frac{A(\beta) -LB(\beta)}{L} >  B(\beta),
\]
that is,  $A(\beta) > 2LB(\beta)$, we have that 
\begin{equation}\label{qnu0}
c_{min}(\beta) \leq \inf_{\nu > \sqrt{B}} \sup_{U \in (0,1)} \left\{ 
\left( \nu - \frac{B(\beta)}{\nu}\right)  h'(U) + \frac{A(\beta)}{\nu}
\right\} = q(\nu_0) =2\sqrt{L}\sqrt{A(\beta)-LB(\beta)}.
\end{equation}

Comparison of $q(\nu_0)$ in (\ref{qnu0}) with $c_{nl}(\beta)$
then shows that $c_{min}(\beta) \neq c_{nl}(\beta)$ if $A(\beta) > 2L
B(\beta)$.
\end{Proof}

Now we can formulate sufficient conditions for minimality
exchange. Below we say that a solution $\beta^*$ of the equation
$A(\beta)=2B(\beta)$ is {\bf non-degenerate} if the graphs of the
functions $A(\cdot)$ and $2B(\cdot)$ intersect transversely at
$\beta^*$.  The following result  applies in all the examples in \cite{CGM,HR}
mentioned above and covers the general case when $h(u)$ is concave and
there is a non-degenerate solution to $A(\beta)=2B(\beta)$.

\begin{theo}\label{me}
Suppose there is a non-degenerate solution $\beta^*$ to the equation
$A(\beta)=2B(\beta)$. Then if
\[
 L= h'(0)=1,
\]
there is a minimality exchange at $\beta=\beta^*$. 
\end{theo}

\begin{Proof}
  Since if $A(\beta) < 2 B(\beta)$ we have that
  $c_{min}(\beta)=c_{nl}(\beta)$ by Lemma \ref{lnl}, and since by
  (\ref{qnu0}) with $L=1$, for all $A(\beta)> 2B(\beta)$,
  $c_{min}(\beta) = c_l(\beta)$, non-degeneracy of the solution
  $\beta^*$ of $A(\beta)=2B(\beta)$ implies that there is an exchange
  of minimality at $\beta^*$.
\end{Proof}

Theorem \ref{me} fully characterises minimality exchange when $L=1$,
that is, when $h'(u)$ attains its supremum $L$ at $u=0$, which holds
in particular when $h$ is concave. If $L > 1$, however, the situation
is less clear. Lemma \ref{lnl} clearly still implies that
$c_{min}(\beta)= c_{nl}(\beta) > c_l(\beta)$, so in particular
nonlinear selection holds, if $A(\beta) < 2 B(\beta)$, and linear
selection holds, with $c_{min}(\beta) = c_{nl}(\beta) = c_l(\beta)$ if
$A(\beta)=2B(\beta)$, but whether it is possible to have again
nonlinear selection for some $\beta$ with $A(\beta) > 2B(\beta)$,
either with the minimal speed corresponding to the explicit solution
or another value, is not obvious.  The estimate (\ref{qnu0}) only
applies when $A(\beta) > 2LB(\beta)$, and even in that range,
(\ref{qnu0}) is no longer sufficient to imply linear selection if
$L>1$.

In Theorem \ref{me2} below, we present a result complementary to
Theorem \ref{me} that makes no assumption on $h$ beyond the hypotheses
in Lemma \ref{solv}, but instead imposes monotonicity conditions on
the dependence of $A$ and $B$ on $\beta$.  This yields a partial
answer to what happens when $L>1$ and $A(\beta) > 2B(\beta)$. We begin
with the following preliminary result, based on \cite[Theorem
2.8]{LMN}, which forms the basis for the alternative sufficient
condition for minimality exchange in Theorem \ref{me2}.

\begin{lemma}
\label{mono}
Suppose that $A(\beta)$ and $B(\beta)$ are each non-decreasing in
$\beta$, and $A(\beta)-B(\beta)$ is non-increasing in $\beta$.  If
$c_{min}(\beta_1) > c_l (\beta_1)$ and $\beta_2> \beta_1$, then
$$c_{min}(\beta_2) > c_l (\beta_2).$$
that is, if nonlinear selection holds for some $\beta_1$, nonlinear
selection also holds for any $\beta_2 > \beta_1$,
\end{lemma}

\begin{Proof}
  We draw on Theorem 2.8 of Lucia, Muratov and Novaga \cite{LMN},
  which says that $c_{min}(\beta) > c_l(\beta)$ if and only if there
  exists $c> c_l(\beta)$ and $u \in H^1_c(\R)$ such that
\begin{equation}
\label{lmn2.8}
\Phi_c^{\beta}[u] := \int_{\R} e^{cx} \left( \frac{1}{2} u_x^2 -
  \int_0^u f(s, \beta)\, ds \right) \, dx  \;\leq \; 0,
\end{equation} 
where $H^1_c(\R)$ is  as defined in the proof of Lemma \ref{lnl}. 

First note that it follows from \cite[Theorem 2.8]{LMN} that since
$c_{min}(\beta_1) > c_l (\beta_1)$, there exists $c > c_l(\beta_1)$
and $u \in H^1_c(\R)$ such that $\Phi_c^{\beta_1}[u] \leq 0$. Then
\begin{eqnarray*}
\Phi_c^{\beta_1}[u] & = & \int_{\R} e^{cx} \left( \frac{1}{2} u_x^2 - \int_0^u f(s, \beta_1) ds \right) \, dx \\ 
& = & \int_{\R}  e^{cx} \left( \frac{1}{2} u_x^2 - \int_0^u h(s) (A(\beta_1)  -B(\beta_1) h'(s)) ds \right) \, dx \\ 
& = & \int_{\R}  e^{cx} \left( \frac{1}{2} u_x^2 - A(\beta_1) \int_0^u h(s) \, ds - \frac{B(\beta_1)}{2} h(u)^2 \right) \, dx\\
& \leq & 0, 
\end{eqnarray*}
as $h(0)=0$, and since $\beta_2> \beta_1$ and $A(\cdot)$ and $B(\cdot)$ are
non-decreasing, we have $A(\beta_2) \geq A(\beta_1)$ and
$B(\beta_2) \geq B(\beta_1)$, so that
$$\Phi_c^{\beta_2}[u] \, \leq \, \Phi_c^{\beta_1}[u] \,\leq \, 0,$$
since $h(s)>0$ for $0<s<1$. Moreover, $A(\cdot)-B(\cdot)$ is
non-increasing, so
$$  c_l(\beta_2) = 2 \sqrt{A(\beta_2) - B(\beta_2)} \, \leq 2 \sqrt{A(\beta_1) - B(\beta_1)} = c_l(\beta_1), $$
and hence
$$c> c_l(\beta_1) \geq c_l(\beta_2).$$
Thus $c> c_l(\beta_2)$ and $\Phi_c^{\beta_2}[u] \leq 0$, and hence \cite[Theorem 2.8]{LMN}  implies that $c_{min}(\beta_2) > c_l (\beta_2)$. 
\end{Proof}

The following is an immediate consequence of Lemma \ref{mono}.

\begin{cor}
\label{mono2}
Suppose that $A(\beta)$ and $B(\beta)$ are each non-decreasing in
$\beta$, and that $A(\beta)-B(\beta)$ is a non-increasing function in
$\beta$.  If $c_{min}(\beta_2) = c_l(\beta_2)$ for some $\beta_2$ and
$\beta_1<\beta_2$, then $c_{min}(\beta_1) = c_l(\beta_1)$.
\end{cor}

We can now prove our second set of sufficient conditions for
minimality exchange.

\begin{theo}
\label{me2}
Suppose that $A(\beta)$ and $B(\beta)$ are each non-decreasing in
$\beta$, and $A(\beta)-B(\beta)$ is non-increasing in $\beta$.  If
there is a non-degenerate solution $\beta^*$ to the equation
$A(\beta)=2B(\beta)$, then there is a minimality exchange at
$\beta=\beta^*$, with $c_{min}(\beta) = c_l(\beta)$ for
$\beta \leq \beta^*$ and $c_{min}(\beta) = c_{nl}(\beta) > c_l(\beta)$
for $\beta > \beta^*$.
\end{theo}

\begin{Proof}
  Note first that $A(\beta)-2B(\beta = [A(\beta) -B(\beta)]-B(\beta)$
  is non-increasing in $\beta$, so since the graphs of $A(\cdot)$ and
  $2B(\cdot)$ intersect transversally at $\beta^*$, it follows that
  $A(\beta) > 2B(\beta)$ when $\beta< \beta^*$, whereas
  $A(\beta) < 2B(\beta)$ when $\beta > \beta^*$.  Lemma \ref{lnl} then
  implies that $c_{min}(\beta)=c_{nl}(\beta)$ when $\beta> \beta^*$,
  whereas Corollary \ref{mono2} implies that linear selection holds
  when $\beta< \beta^*$.
\end{Proof}

Note that for the two concrete examples of minimality exchange
discussed in Section 3, both Theorem \ref{me} and Theorem \ref{me2}
apply.

An example of a solvable problem for which Theorem \ref{me2} applies
but Theorem \ref{me} does not, is given by taking $A=1$, $B=\beta/2$
and $h(u)=e^{2u}u(1-u)$, which is not concave. Then $L=1.52218$,
$c_l= \sqrt{4-2\beta}$, $c_{nl}=\sqrt{2/\beta}$,
$c_l(\beta) =c_{nl}(\beta)$ at $\beta^*=1$, and Theorem \ref{me2}
ensures that there is minimality exchange at $\beta^*=1$.

\section{Conclusions}
In this article we have focussed on a class of parameter-dependent
monostable reaction-diffusion equations with explicit travelling-wave
solutions and used this class to explore the phenomenon of minimality
exchange, when the minimal wave speed switches from a linearly
determined value to the speed of the explicitly determined front as a
parameter changes. Two alternative sets of sufficient conditions for
minimality exchange are proved, in Theorems \ref{me} and
\ref{me2}. Why there should be such an exchange, not only from linear
selection to nonlinear selection, but to nonlinear selection given by
an {\em explicit} solution, is quite puzzling at first sight.  Our
framework here provides insight into why minimality exchange of this
type occurs, and includes concrete examples from \cite{CGM, GK, HR,
  LMN}. The proofs owe much to a variety of tools for determining
whether there is linear or nonlinear selection - in particular, ideas
developed previously in the special case of an isotropic
liquid-crystal model \cite{CGM}, as well as general results from
\cite{HR, LMN}. Some additional interesting material about minimal
wave speeds is given in \cite[Section 10.1.1]{GK}, including Theorem
10.12, which provides sufficient criteria that can be used to identify
cases when a given explicit solution has the minimal wave speed, and
the examples that follow.

We have treated one class of parameter-dependent solvable equations
that includes important special cases, but clearly there are many
further solvability results for explicit travelling-wave solutions in
the literature. See, for instance, \cite[Chapter 13]{GK} and
\cite{PL}. In addition, the change of variables $G:=1/F$ converts
(\ref{feq}) into an Abel equation, for which certain classes of
explicit solutions can be found using tools such as the Chiellini
integrability condition and the Lemke transformation (see, for
example, \cite{HM} and the references therein). It would be
interesting to expand and develop the approach introduced here to
cover a larger range of explicit solutions to obtain further insight
into the mechanisms for minimality exchange.

\end{document}